\documentclass[%
 aip,
cp,  
 amsmath,amssymb,
 reprint,%
]{revtex4-2}

\usepackage{graphicx}
\usepackage{dcolumn}
\usepackage{bm}

\usepackage[utf8]{inputenc}
\usepackage[T1]{fontenc}
\usepackage{mathptmx} 

\begin{document}

\title{Research of the hereditary dynamic Riccati system with modification fractional differential operator of Gerasimov-Caputo}

\author{Tverdyi D.A.} 
 \email{dimsolid95@gmail.com}
 \affiliation{Institute of Applied Mathematics and Automation KBSC RAS, Nalchik, Russia;}
 \affiliation{KamGU Vitus Bering, Petropavlovsk-Kamchatsky, Russia;}
\author{Parovik R.I.}%
 \email{romanparovik@gmail.com}
  \affiliation{Institute of Cosmophysical Research and Radio Wave Propagation FEB RAS, Kamchatka Territory, Paratunka village, Russia;}

\date{\today} 

\begin{abstract}
In this paper, we study the Cauchy problem for the Riccati differential equation with constant coefficients and a modified Gerasimov-Caputo type fractional differential operator of variable order. Using Newton's numerical algorithm, calculation curves are constructed taking into account different values of the Cauchy problem parameters. The calculation results are compared with the previously obtained results. The computational accuracy of the numerical algorithm is investigated. It is shown using the Runge rule that the computational accuracy tends to the accuracy of the numerical method when increasing the nodes of the calculated grid. 
\end{abstract}

\maketitle

\section{\label{sec:level1}Introduction}

The classical Riccati equation describes processes with saturation, for example, the logistic law in biology or in Economics, and the Riccati equation occurs after some transformations in physical problems, for example, in problems of reflection of waves from a non-uniform surface. The Riccati equation is well studied and can be found in various reference books and encyclopedias on mathematics or in the book \cite{1ReidWT}.

However, in the last decade there has been a surge in the application of fractional calculus to the study of the Riccati equation \cite{1ReidWT,2KhashanMM,3HouJ,4SakarMG,5MerdanM,6KhaderMM,7KhaderMM,8GoharM,9SyamMI,10EzzEldienSS,11SalehiY,12KhanNA,13EzzEldienSS,14SweilamNH,15AminikhahH,16SyamM,17TverdyiDA,18TverdyiDA,19TverdyiDA,20TverdyiDA,21TverdyiDA}. This is due to the fact that the processes with saturation may have effects hereditarily. The effect of hereditarily means that the system or process can remember about their background and from the point of view of mathematics can be described by means of integro-differential equations with difference kernels with memory functions \cite{22VolterraV}. When choosing power functions of memory, we naturally proceed to the well-known mathematical apparatus of fractional calculus, in particular to derivatives of fractional orders \cite{23KilbasAA,24OldhamKB,25MillerKS}.  The Riccati equation with a fractional order derivative is called the Riccati fractional equation.

The fractional Riccati equation with constant order was studied in \cite{2KhashanMM,3HouJ,4SakarMG,5MerdanM,6KhaderMM,7KhaderMM,8GoharM,9SyamMI,10EzzEldienSS,11SalehiY,12KhanNA,13EzzEldienSS,14SweilamNH,15AminikhahH}. However, a broader class of fractional Riccati equations with variable order is of interest \cite{16SyamM,17TverdyiDA,18TverdyiDA,19TverdyiDA,20TverdyiDA,21TverdyiDA}.

In the works of the authors \cite{17TverdyiDA,18TverdyiDA,19TverdyiDA,20TverdyiDA,21TverdyiDA}, a fractional Riccati equation of variable order with non-constant coefficients was studied. Using the Newton-Raphson numerical method, calculated curves were obtained and the accuracy of the numerical method was investigated. The research results were used for modeling some logistic laws in the model of solar activity dynamics \cite{21TverdyiDA}. Calculations and visualizations in these works were performed in the Maple computer mathematics environment.

In this paper, we will use a modified Gerasimov-Caputo type fractional derivative operator of variable order. Next, we will study the numerical solution of the fractional Riccati equation by analogy with the work \cite{20TverdyiDA} and compare the results obtained.

\section{Some basic definitions}
Here we will consider the main definitions from the theory of fractional calculus, and its aspects can be studied in more detail in books \cite{23KilbasAA,24OldhamKB,25MillerKS}.

\textbf{Definition 1:} The Gerasimov-Caputo fractional derivative of variable order $\alpha(t)$ has the form:
\begin{eqnarray}
\partial^{\alpha(t)}_{0t}u(\tau)=
\left\{
\begin{array}{ll}
\frac{1}{\Gamma(m-\alpha(t))} \int_{0}^{t} \frac{x^{(m)}(\tau)d\tau}{(t-\tau)^{\alpha(t)+1-m}},\qquad 0\leq m-1 < \alpha(t) < m\\
\frac{d^{m}x(t)}{dt^{m}}, \qquad m\in N.
\end{array}
\right.
\label{eq:1}
\end{eqnarray}
where $\Gamma(z)=\int_{0}^{\infty}t^{z-1}\exp^{-t}dt, Re(z)>0$ -- Euler's gamma function.

\textbf{Definition 2:} Modified Gerasimov-Caputo fractional derivative of variable order $\gamma(t-\tau)$ has the form:
\begin{eqnarray}
\partial^{\gamma(t-\tau)}_{0t}u(\tau)=
\left\{
\begin{array}{ll}
\frac{1}{\Gamma(m-\gamma(t-\tau))} \int_{0}^{t} \frac{x^{(m)}(\tau)d\tau}{(t-\tau)^{\gamma(t-\tau)+1-m}},\qquad 0\leq m-1 < \gamma(t-\tau) < m\\
\frac{d^{m}x(t)}{dt^{m}}, \qquad m\in N.
\end{array}
\right.
\label{eq:2}
\end{eqnarray}
Derivatives of fractional variables of orders (1) and (2) are given in \cite{26LorenzoCF}, and some of their properties are considered there.

\textbf{Note 1:} Note that if the orders of fractional derivatives in (\ref{eq:1}) and (\ref{eq:2}) are constants, then they coincide and the operators are Gerasimov-Caputo derivatives of constant orders \cite{27GerasimovA,28CaputoM}.

\section{Task definition}
Consider the following hereditary equation that is the analogue of the equation of  Rikkati:
\begin{eqnarray}
\int_{0}^{t}K(t-\tau) \dot{u}(\tau)d\tau+a(t)u^{2}(t)+b(t)u(t)+c(t)=0
\label{eq:3}
\end{eqnarray}
where $u(t)\in C[0,T]$ -- function of a solution, $K(t-\tau)$ -- memory function, $t\in[0,T]$ -- current time, $T>0$ -- simulation time, $a(t), b(t), c(t)$ -- coefficients set by the function.

\textbf{Note 2:} Note that if the memory function $K(t-\tau)$ is a Heaviside function, then the process has full memory; if it is a Dirac function, then there is no memory. 

We will consider an intermediate case when the system gradually “forgets” its background over time. To do this, select the memory function as follows:
\begin{eqnarray}
K(t-\tau) = \frac{(t-\tau)^{-\alpha(t)}}{\Gamma(1-\alpha(t))}, \qquad 0<\alpha(t)<1
\label{eq:4}
\end{eqnarray}
where $\alpha(t)$ -- a function that is responsible for the intensity of the process under study, then taking into account the Definition of 1 for $m=1$, we come to the following fractional Riccati equation of variable order:
\begin{eqnarray}
\partial^{\alpha(t)}_{0t}u(\tau)+a(t)u^{2}(t)+b(t)u(t)+c(t)=0
\label{eq:5}
\end{eqnarray}
for which the initial local condition is valid:
\begin{eqnarray}
u(0)=u_{0},
\label{eq:6}
\end{eqnarray}
Equation (\ref{eq:5}) and initial condition (\ref{eq:6}) form the Cauchy problem for the fractional Riccati equation of variable order. It was studied by the authors in \cite{17TverdyiDA}.

Let us introduce into consideration a different function of memory:
\begin{eqnarray}
K(t-\tau) = \frac{(t-\tau)^{-\gamma(t-\tau)}}{\Gamma(1-\gamma(t-\tau))}, \qquad 0<\gamma(t-\tau)<1
\label{eq:7}
\end{eqnarray}
Here, the order of the $\gamma(t-\tau)$ function with a lagging argument also determines the intensity of the process under consideration.

Taking into account (\ref{eq:7}) and Definition 2 for $m=1$, the fractional Riccati equation takes the form:
\begin{eqnarray}
\partial^{\gamma(t-\tau)}_{0t}u(\tau)+a(t)u^{2}(t)+b(t)u(t)+c(t)=0
\label{eq:8}
\end{eqnarray}
In the future, we will investigate the Cauchy problem (\ref{eq:8}), (\ref{eq:6}).

\section{The method of solving}
Since the Cauchy problem (\ref{eq:8}) and (\ref{eq:6}) in General does not have an exact solution, we will use numerical methods to solve it. To do this, we divide the time interval $[0,T]$ into $N$ equal parts (grid nodes), where $h=\frac{T}{N}$ is the sampling step $t_{n}=nh, n=0,\dots,N-1$, and the solution function $u(t_{n})=u_{n}$.

The approximation of the fractional derivative in equation (8) will take the form:
\begin{eqnarray}
\partial^{\gamma(t-\tau)}_{0t} \approx \sum\limits_{i=1}^{k} \omega_{i,\gamma} \left(u_{k-i+1}-u_{k-i}\right), \qquad k=1,\dots,N
\label{eq:9}
\end{eqnarray}
where,  $\omega_{i,\gamma}$ weight coefficients of the quadrature formula of trapezes.

Note that the method of approximating the Cauchy problem (\ref{eq:5}), (\ref{eq:6}) is based on the results of \cite{29ParovikRI}.

Based on (\ref{eq:9}), the Cauchy problem (\ref{eq:8}) and (\ref{eq:6}) can be rewritten in the difference statement:
\begin{eqnarray}
 \sum\limits_{i=1}^{k} \omega_{i,\gamma} \left(u_{k-i+1}-u_{k-i} \right)+a_k u_k^2+b_k u_k+c_k=0, \qquad k=1,\dots,N
\label{eq:10}
\end{eqnarray}
For the solution, the iterative Newton-Raphson method was used, which, as a rule, provides fast convergence.
The algorithm consists of the following steps:
\begin{enumerate}
\item Define a function for forming a Jacobian:
\begin{eqnarray*}
f_{k} = \sum\limits_{i=1}^{k}\omega_{i,\gamma}\left(u_{k-i+1}-u_{k-i} \right) + a_{k} u^{2}_{k} + b_{k}u_{k} + c_{k}, \qquad k=1,\dots,N
\end{eqnarray*}
\item Define the elements of the Jacobian:
\begin{eqnarray*}
R_{n,m} = \frac{d f_{n}}{d u_{m}}, \qquad n=1,\dots,N, m=1,\dots,N
\end{eqnarray*}
\item Let's write a matrix iterative equation of the form:
\begin{eqnarray}
U_{m+1} = U_{m} - \frac{1}{J(U_{m})} F(U_{m})
\label{eq:11}
\end{eqnarray}
where
\begin{eqnarray*}
U_m=\left( \begin{array}{c} u_1 \\ u_2 \\ \vdots \\ u_N \end{array} \right),
F(U_m)=\left( \begin{array}{c} f_1 \\ f_2 \\ \vdots \\ f_N \end{array} \right),
J(U_m)=\left( \begin{array}{cccc} R_{1,1} & R_{1,2} & \ldots & R_{1,N} \\ R_{2,1} & \ldots & \ldots & R_{2,N} \\ \vdots & \vdots & R_{3,3} & \vdots \\ R_{N,1} & \ldots & \ldots & R_{N,N} \end{array} \right)
\end{eqnarray*}
\item We start the iterative process until $r>\epsilon$, where $\epsilon = 10^{-4}$ -- accuracy, $r = 10^{3} \epsilon $ -- stop criterion, and $\| J(U_m) \| \neq 0$ -- convergence criterion of the method, where $r = \| U_{m+1}-U_{m} \| \leq \epsilon$ -- stability criterion.
\end{enumerate}

\textbf{Note 3:} It should be noted that the main difficulty in applying this method is to find the inverse Jacobi matrix in (\ref{eq:11}). To solve this problem, we can use the Gauss – Jordan method \cite{30LipschutzS}.

\section{The simulation parameters}
The variable order $\gamma(t-\tau)$ is considered as a periodic function:
\begin{eqnarray}
\gamma(t-\tau) = \frac{\theta cos(\mu h (t-\tau)) + 2\delta}{2}
\label{eq:12}
\end{eqnarray}
where $\delta$ -- the shift coefficient that allows the condition to be met  $0<\gamma(t-\tau)<1$, $\theta$ -- oscillation amplitude, and $\mu$ -- oscillation frequency.

The coefficients in the difference equation (\ref{eq:10}) will take the following values:
\begin{eqnarray}
a_{k} = -\frac{k}{N}, b_{k} = 0, c_{k} =\frac{k}{N}
\label{eq:13}
\end{eqnarray}

The introduction of coefficients of the form (\ref{eq:13}), as shown in the author's work \cite{21TverdyiDA}, is due to the appearance of distribution curves similar to the S-shaped (logistics) curve, which has its own application.

Consider the following examples for $T = 50, N = 2000$:
\begin{enumerate}
\item Example: $\alpha = \gamma = const = 1 \approx 0.9999$;
\item Example: $0.5 < \alpha,\gamma < 1$ with parameters (\ref{eq:12}): $\delta=0.75, \theta=0.5, \mu=\frac{\pi}{2}$;
\item Example: $0 < \alpha,\gamma < 1$ with parameters (\ref{eq:12}): $\delta=0.5, \theta=0.5, \mu=\frac{\pi}{2}$;
\item Example: $0 < \alpha,\gamma < 0.5$ with parameters (\ref{eq:12}): $\delta=0.25, \theta=0.5, \mu=\frac{\pi}{2}$;
\end{enumerate}
We study the solution of Cauchy problems (\ref{eq:5}), (\ref{eq:6}) and compare it with the solution of Cauchy problem (\ref{eq:8}), (\ref{eq:6}).

\section{Simulation result}
We perform verification and Show that using fractional derivatives (\ref{eq:1}) of the form $\partial^{\alpha(t)}_{0t}u(\tau)$ and (\ref{eq:2}) of the form $\partial^{\gamma(t-\tau)}_{0t}u(\tau)$, for $\alpha = \gamma = const = 1$, we get the same distribution curves.
\begin{figure}[h!]
\includegraphics[scale=0.5]{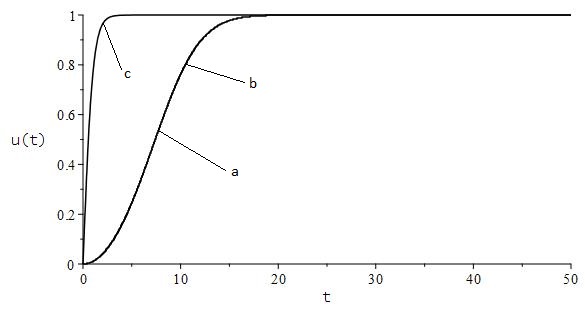}
\caption{\label{fig:1} Verification. \textbf{a,b)} Example 1 for the operator (\ref{eq:1}) and (\ref{eq:2}), \textbf{c)} The classic solution, i.e. when $a_{k} = -1, b_{k} = 0, c_{k} =1$, and $\alpha = \gamma = const = 1$.}
\end{figure}
\begin{figure}[h!]
\includegraphics[scale=0.5]{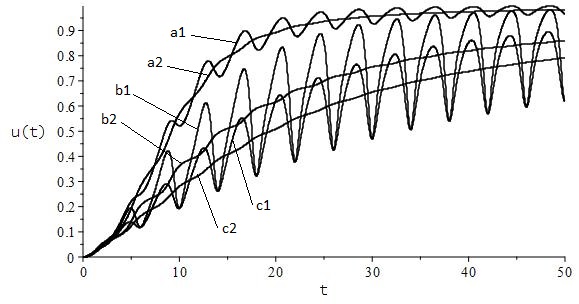}
\caption{\label{fig:2} Distribution curves: \textbf{a1)} Example 2 for (\ref{eq:1}), \textbf{a2)} Example 2 for (\ref{eq:2}), \textbf{b1)} Example 3 for (\ref{eq:1}), \textbf{b2)} Example 3 for (\ref{eq:2}), \textbf{c1)} Example 4 for (\ref{eq:1}), \textbf{c2)} Example 4 for (\ref{eq:2}).}
\end{figure}

Consider (in table \ref{tab:table1}, \ref{tab:table2}, \ref{tab:table3}, \ref{tab:table4}) the change in the absolute error $\epsilon$ and the calculated order of accuracy $p$ for the numerical Newton-Rapson algorithm for solutions of Cauchy problems (\ref{eq:5}), (\ref{eq:6}) and (\ref{eq:8}), (\ref{eq:6}) when the sampling step $h$ is reduced. To calculate the absolute error $\epsilon$, we will use the Runge rule \cite{32BerezinIS}:
\begin{eqnarray}
\epsilon = max \left( \frac{\arrowvert u_{2N,2k-1}-u_{N,k} \arrowvert}{2^{p_{aprior}-1}} \right), \qquad k=1,\dots,N, p_{aprior}=1
\label{eq:14}
\end{eqnarray}

\textbf{Note 4:} A priori accuracy $p_{aprior}$ of the solution in this method is set to 1. This follows from the General scheme approximation order given at the boundary nodes of the grid.

The accuracy of the solution was calculated using the following formula:
\begin{eqnarray}
p = \log_{\frac{h_1}{h_2}} \left( \frac{\epsilon_{1}}{\epsilon_{2}} \right),
\label{eq:15}
\end{eqnarray}
where $h_1,h_2 = \dfrac{h_1}{2}$ -- steps of sampling, $\epsilon_{1},\epsilon_{2}$ -- errors on steps $h_1$ (previous) and $h_2$ (current).

\begin{table}[h!]
\caption{\label{tab:table1} Example 1: $\alpha = \gamma = const = 1$}
\begin{ruledtabular}
\begin{tabular}{cccccc}
 &\multicolumn{1}{c}{$T=50$}&\multicolumn{2}{c}{$\alpha(t)$}&\multicolumn{2}{c}{$\gamma(t-\tau)$}\\
 $N$& $h$ & $\epsilon$ & $p$ & $\epsilon$ & $p$\\ \hline
 129 & 0.387 & 0.063871 & -        & 0.063871 & - \\
 259 & 0.193 & 0.032515 & 0.974045 & 0.032515 & 0.974045 \\
 519 & 0.096 & 0.016398 & 0.987562 & 0.016398 & 0.987562 \\
 1039& 0.048 & 0.008233 & 0.993947 & 0.008233 & 0.993947 \\
 2079& 0.024 & 0.004125 & 0.997027 & 0.004125 & 0.997027 \\
\end{tabular}
\end{ruledtabular}
\end{table}

\begin{table}[h!]
\caption{\label{tab:table2} Example 2: $0.5 < \alpha,\gamma < 1$ with parameters (\ref{eq:12}): $\delta=0.75, \theta=0.5, \mu=\frac{\pi}{2}$}
\begin{ruledtabular}
\begin{tabular}{cccccc}
 &\multicolumn{1}{c}{$T=50$}&\multicolumn{2}{c}{$\alpha(t)$}&\multicolumn{2}{c}{$\gamma(t-\tau)$}\\
 $N$& $h$ & $\epsilon$ & $p$ & $\epsilon$ & $p$\\ \hline
 129 & 0.387 & 0.070173 & -        & 0.045638 & - \\
 259 & 0.193 & 0.034098 & 1.041212 & 0.023454 & 0.960369 \\
 519 & 0.096 & 0.017016 & 1.002759 & 0.011735 & 0.998989 \\
 1039& 0.048 & 0.008363 & 1.024701 & 0.005831 & 1.008889 \\
 2079& 0.024 & 0.004117 & 1.022315 & 0.002892 & 1.011413 \\
\end{tabular}
\end{ruledtabular}
\end{table}

\begin{table}[h!]
\caption{\label{tab:table3} Example 3: $0 < \alpha,\gamma < 1$ with parameters (\ref{eq:12}): $\delta=0.5, \theta=0.5, \mu=\frac{\pi}{2}$}
\begin{ruledtabular}
\begin{tabular}{cccccc}
 &\multicolumn{1}{c}{$T=50$}&\multicolumn{2}{c}{$\alpha(t)$}&\multicolumn{2}{c}{$\gamma(t-\tau)$}\\
 $N$& $h$ & $\epsilon$ & $p$ & $\epsilon$ & $p$\\ \hline
 129 & 0.387 & 0.305098 & -        & 0.028593 & - \\
 259 & 0.193 & 0.182349 & 0.742568 & 0.013625 & 1.069381 \\
 519 & 0.096 & 0.095837 & 0.928038 & 0.006677 & 1.029022 \\
 1039& 0.048 & 0.048632 & 0.978686 & 0.003289 & 1.021300 \\
 2079& 0.024 & 0.024420 & 0.993831 & 0.001633 & 1.009579 \\
\end{tabular}
\end{ruledtabular}
\end{table}

\begin{table}[h!]
\caption{\label{tab:table4} Example 4: $0 < \alpha,\gamma < 0.5$ with parameters (\ref{eq:12}): $\delta=0.25, \theta=0.5, \mu=\frac{\pi}{2}$}
\begin{ruledtabular}
\begin{tabular}{cccccc}
 &\multicolumn{1}{c}{$T=50$}&\multicolumn{2}{c}{$\alpha(t)$}&\multicolumn{2}{c}{$\gamma(t-\tau)$}\\
 $N$& $h$ & $\epsilon$ & $p$ & $\epsilon$ & $p$\\ \hline
 129 & 0.387 & 0.180735 & -        & 0.016893 & - \\
 259 & 0.193 & 0.110282 & 0.712672 & 0.008336 & 1.019414 \\
 519 & 0.096 & 0.056046 & 0.976513 & 0.004205 & 0.986786 \\
 1039& 0.048 & 0.028572 & 0.972018 & 0.002139 & 0.975103 \\
 2079& 0.024 & 0.014311 & 0.997402 & 0.001086 & 0.976774 \\
\end{tabular}
\end{ruledtabular}
\end{table}

\newpage

\section{Conclusion}
The Cauchy problem for the fractional Riccati equation (\ref{eq:8}), (\ref{eq:6}) with variable coefficients was considered. The proposed Cauchy problem is numerically analyzed using the Newton-Raphson method. The computational accuracy of the Runge rule for the applied numerical method is estimated. It is shown that the computational accuracy tends to the accuracy of the method as the number of calculated nodes increases. Next, we compared the solution of the Cauchy problem (\ref{eq:5}), (\ref{eq:6}) with the Cauchy problem (\ref{eq:8}), (\ref{eq:6}), which differ from each other, but maintain a common trend. This fact suggests that the Cauchy problem (\ref{eq:8}), (\ref{eq:6}) takes place and its solution can be used in the study of processes with saturation and memory effects, for example, in forecasting economic cycles and crises \cite{31MakarovDV}.

\bibliography{aipsamp}

\end{document}